\input amstex
\documentstyle{amsppt}
\magnification=1200
\hsize=13.8cm
\vsize=19.0 cm
\catcode`\@=11
\def\NoLogo{\let\logo@\empty}
\catcode`\@=\active
\NoLogo

\def\e{\epsilon}

\def \b {\beta}

\def\lf{\left}
\def\ri{\right}
\def\a{\alpha}
\def\g{\gamma}
\def\p{\partial}
\def\delbar{\bar\delta}
\def\ddbar{\partial\bar\partial}
\def\dbar{\bar\partial}

\def\l{\lambda}
\def\L{\Lambda}
\def\C{\Bbb C}
\def\R{\Bbb R}

\def\vp{\varphi}

\def\dbar{\bar\partial}

\def\bb{{\bar\beta}}

\def\abb{{\alpha\bar\beta}}

\def \D {\Delta}
\def\heat{\lf(\frac{\p}{\p t}-\Delta\ri)}
\def\aint{\frac{\ \ }{\ \ }{\hskip -0.4cm}\int}

\def\bbar{{\bar \beta}}

\leftheadtext{Lei Ni }
\rightheadtext{Monotonicity and K\"ahler-Ricci flow}

\topmatter

\title{Monotonicity and K\"ahler-Ricci flow}\endtitle

\author{Lei Ni\footnotemark}\endauthor
\footnotetext{Research partially supported by NSF grant DMS-0203023.}

\address
Department of Mathematics, University of California, la Jolla, CA 92093
\endaddress
\email{
lni\@math.ucsd.edu}
\endemail

\affil
{
University of California, San Diego
}
\endaffil

\date  October, 2002
\enddate

\endtopmatter

\document

\subheading{\S0 Introduction}

In this paper, we shall give a geometric account of the linear trace
Li-Yau-Hamilton (which will be abbreviated as LYH)
inequality for the K\"ahler-Ricci flow
proved by Luen-Fai Tam and the author in [NT1].
To put the result,  especially the Liouville theorem for the plurisubharmonic
functions, into the right perspective we would also describe  some
dualities existed in both linear and nonlinear analysis first.
The purpose is
to show the reader that the linear trace LYH for the K\"ahler-Ricci flow
can be thought as a `global' parabolic version of the classical
monotonicity formula for the analytic hypersurfaces in $\C^m$ (or more
generally for the positive $(1,1)$ currents). We also include some new
results. For example the Harnack inequality for the nondivergent
elliptic operator
on complete K\"ahler manifolds with nonnegative holomorphic bisectional
curvature was proved in  Theorem 1.5. Another is the
LYH inequality for the Hermitian-Einstein flow
coupled with K\"ahler-Ricci flow (cf. Theorem 3.5).  We also prove that the
sufficient and necessary condition for the equality in the linear trace
LYH inequality is that the solution is a K\"ahler-Ricci soliton. (It has
no restriction on the tensor satisfying the linear Lichnerowicz-Laplacian
heat equation.) See, Theorem 4.1 and Theorem 4.2. A direct consequence of
these result is that type II (III) limit solutions of Ricci (K\"ahler-Ricci)
flow are gradient
 (expanding)  solitons (K\"ahler-Ricci solitons). Theorem 4.1 and
 Theorem 4.2
 generalizes  and unifies the previous theorems (cf. [H3],
[Ca2], [C-Z]) of Hamilton, Cao and more recently Chen-Zhu on the
limit solutions to the Ricci flow considerably.

There are monotonicity formulae for the parabolic equations such as the
harmonic map heat equation and mean
curvature flow (cf. [Hu], [St] and [H2]).
But in the author's point of view
they  all are  still `local' in the sense that the precise monotonicity
 only holds for (or inside)
locally symmetric manifolds. Another important distinction is that
 the `global monotonicity' (which holds on complete K\"ahler manifolds with
 nonnegative bisectional curvature)
 derived from
the LYH inequality here
is a  point-wise estimate instead of the monotonicity of an
integral quantity as in the
previous mentioned cases such as those in [St] and [Hu].
There have been some other
geometric interpretations on LYH inequality proved by Hamilton, for example
 in [C-C1],
as well as  on the  Chow-Hamilton's linear trace LYH  in [C-C2].
This  observation here
relates the minimal submanifold theory, in particular the
monotonicity of the volume,  to the LYH inequality for Ricci flow, at least
for complex analytic case, for the first time. Hopefully this will introduce
further work along this direction and bring out
more understanding and results
to the study of Ricci flow. Let us start with a duality for linear elliptic
PDE,
which the author learned from Professor Peter Li during his graduate
study at UCI back in 1997-1998 (cf. [L-W] for more works in this direction).
This duality is hinged on the Harnack
inequality.

\medskip

\subheading{\S1 Harnack inequality}

On $\R^n$, the celebrated De Giorgi-Nash Moser theory proves the
Harnack inequality for the uniformly elliptic operator of divergence form.

\proclaim{Theorem 1.1}
Let $L=\sum_{ij}\p_i (a_{ij}(x)\p_j) $ be a uniformly elliptic operator on
$\R^n$ with
$\l I\le \lf(a_{ij}(x)\ri)\le \L I$. Here constants $\L\ge \l>0$. Let
$u\ge 0$ be a (local) $W^{1,2}$-weak solution to  $Lu=0$ in $B_o(2R)$.
Then there exists a constant
$C=C(n, \frac{\L}{\l})>0$ such that
$$
 \sup_{B_o(R)}u \le C\inf_{B_o(R)}u. \tag 1.1
$$
\endproclaim
The duality provided by the Harnack inequality is between the local
$C^\alpha$-regularity of any $W^{1,2}$-solution and the Liouville theorem
for the positive solution of
$Lu=0$. As it is well-known (cf. [G-T]) that applying (1.1) to balls with
small radius $R$ implies that any $W^{1,2}$ solution to $Lu=0$ is $C^{\a}$.
It is also easy to see that if $u\ge 0$ is a solution of $Lu=0$ on $\R^n$
then $u$ must be a constant. This can be seen by applying (1.1) to
$u-\inf_{\R^n}u$ on $B_o(R)$ and letting $R\to \infty$.

When study the analysis on general manifolds, since any manifold locally
looks Euclidean and the metric is quasi-isometric to the Euclidean one the
$C^\alpha$-regularity still holds for any week solution. One can also have
the Harnack. However, the Harnack estimating constant
$C=C(n,\frac{\L}{\l})$ depends on the local geometry of $M$, namely the
metric tensor $g_{ij}$, more precisely the quotient of
the maximum eigenvalue over the minimal eigenvalue,
on local coordinate charts. This is a serious
drawbacks since  in order to study
the  nonlinear problems one also wish to obtain the Harnack inequality
 with the
estimating constant depending only on some global geometric-analytic
 information such as the lower bound of the curvature or
the  the behavior of the volume functional, etc.
In a certain sense, one should think
of the metric $g_{ij}$ as being unknown in details, but satisfying
certain geometric (coordinate invariant) conditions.
The corresponding theory on general manifolds starts with the work of
Bombieri and Giusti [B-G], where they carried out the geometric analysis of
 Moser's proofs to minimizing submanifolds and proved the Harnack inequality
for the uniformly elliptic operator of divergence form on such manifolds.
Later in [Y], a gradient estimate method was initiated and a similar Harnack
inequality as in Theorem 1.1 was derived out of a gradient
estimate for positive harmonic functions on manifolds with nonnegative Ricci
curvature. In particular, in [C-Y], the following theorem on
harmonic functions was proved.

\proclaim{Theorem 1.2}
Let $M$ be a complete Riemannian manifold with nonnegative Ricci curvature.
Let $u$ be a harmonic function on $B_o(2R)$. Then there exists a constant
$C=C(n)>0$ such that
$$
\sup_{B_o(R)}|\nabla u|^2\le \frac{C(n)}{R^2}\sup_{B_o(2R)}u^2. \tag 1.2
$$
In particular, any global harmonic function (defined on $M$)
of the sublinear growth is a constant.
\endproclaim

The result was later also established for general uniformly
elliptic operator of divergence form in [Gr] and [Sa] independently.
In fact, they proved the Harnack for the corresponding parabolic equations
under some  geometric assumptions which are than the nonnegativity
 of the Ricci
curvature. For the sake of this paper we just state their result
in the following weaker form for the elliptic operators.

\proclaim{Theorem 1.3}
Let $M$ be a complete Riemannian manifold with nonnegative Ricci curvature.
 Let $L$ be a uniformly elliptic operator of divergence form. Let $u\ge 0$
be a solution to $Lu=0$ on $B_o(2R)$. Then there exists a constant
$C=C(n, L)$ such that
$$
\sup_{B_o(R)}u\le C\inf_{B_o(R)}u. \tag 1.3
$$
In particular, any positive global $L$-harmonic function is a constant.
\endproclaim

There are also many developments in establishing the Harnack inequality for
the uniform elliptic operator of nondivergence form. The theory for the
Euclidean (local) case was developed by Krylov and Safonov based on the
Alexanderoff maximum principle (cf. [G-T], Chapter 8 for the detailed
theory).
We just want  to mention a corresponding
global
result for the manifolds with nonnegative Riemannian sectional curvature
which was proved by Cabr\'e [C] more recently.

\proclaim{Theorem 1.4}
Let $M$ be a complete Riemannian manifold with nonnegative sectional
curvature. Let $L$ be a uniformly elliptic operator of nondivergence form.
Then (1.3) still holds for nonnegative solution $u$ on $B_o(2R)$.
\endproclaim

Using a comparison theorem  established in [C-N], we shall
 extend
the above result to the complete K\"ahler  manifolds with nonnegative
holomorphic bisectional curvature. More precisely, let $L$ be a
elliptic operator which locally can be written as $L=\sum_{\abb}a^{\abb}(z)
\frac{\p^2}{\p z_{\a}\p \overline{z}_{\beta}}$ such that there exists
$\L\ge \l >0$  with
$\L |\eta|^2\ge a^{\abb}(z)\eta_{\a}
\overline{\eta_{\beta}}\ge \l |\eta|^2$ for
any $\eta=\eta_{\a}dz^{\a}\in T^*_z M$ over any $z\in M$. The following
result can be viewed as a generalization of Theorem 1.4.

\proclaim{Theorem 1.5}
Let $M$ be  a complete K\"ahler manifold (of complex
dimension m) with nonnegative holomorphic bisectional curvature.
Let $u$ be a smooth function in $B_o(2R)$ satisfying $u\ge 0$ in $B_o(2R)$.
Then
$$
\sup_{B_o(R)}u\le C\lf\{ \inf_{B_o(R)}u+\frac{R^2}{V_o^{\frac{1}{n}}(2R)}
\|Lu\|_{L^n(B_o(2R))}\ri\} \tag 1.6
$$
where $C=C(m, \l,\L)>0$. In particular, (1.3) holds for nonnegative solution
to $Lu=0$ and any positive solution $u$ of $Lu=0$ is a constant.
\endproclaim

As in [C] the key is to prove the following proposition.

\proclaim{Proposition 1.1}
Let $u$ be a smooth function in a ball $B_o(7R)$ satisfying $u\ge 0$ in
$B_o(7R)\setminus B_o(5R)$ and $\inf_{B_o(2R)}u\le 1$.
Then
$$
V_o(R)\le \frac{(32)^{m}}{\l^{2m}}\int_{\{u\le 6\}\cap B_o(5R)}\lf\{\lf(
\frac{1}{2m}R^2Lu+\frac{\L}{4}\ri)^{+}\ri\}^{2m}. \tag 1.7
$$
\endproclaim
Once the above result holds one can use the argument of
 section 5-7 in [C], which only uses the nonnegativity of Ricci
curvature, to obtain Theorem 1.5. To prove the above proposition we
need the following two results.

\proclaim{Proposition 1.2}
For any $y\in M$, let $d_y(x)$ be the distance function
from $x$ to $y$. Then
$$
(d^2_y)_{\abb}\le g_{\abb} \tag 1.8
$$
in the distribution sense. In particular, if $x$ does not lie
inside the cut locus of $y$ the above holds point-wisely. Here
$g_{\abb}$ is the K\"ahler metric tensor.
\endproclaim
See [C-N] for the  proof of this comparison result. The following
lemma is from [C].

\proclaim{Lemma 1.1} Let $M$ be a Riemannian manifold.
Let $v$ be a smooth function on $\Omega\subset M$. Consider the map
$\phi: \Omega \to M$ defined by
$$
\phi(z)=\text{exp}_{z}\nabla v(z).
$$
Let $x\in \Omega$ and suppose that $\nabla v\in U_x$. Here
$U_x =\{t\eta:\eta\in T_xM, |\eta|=1, 0\le t\le c(\eta)\}$, where
$c(\eta)=
\sup\{t>0:\text{exp}_xs\eta \text{ is minimizing geodesic on } [0, t]\}
\le \infty$. Set $y=\phi(x)$. Then
$$
\text{Jac}(\phi(x))=\text{Jac}( \text{exp}_x(\nabla v(x))\cdot |\det (D^2(v+
\frac{d_y^2}{2}))(x)|,
$$
where $\text{Jac} (\text{exp}_{x}(\nabla v(x)))$ denotes the Jacobian of
the exponential map at the point $\nabla v(x)\in T_xM$.
\endproclaim
This is just the Lemma 3.2 from [C]. The rest will be devoted to the
proof of Proposition 1.1. The argument is just  an adaptation of  the one
given in [C] to the case with only control on the complex Hessian of
the distance function. The key is the following linear algebra lemma

\proclaim{Lemma 1.2} Let $A$ be a $2m\times 2m$
positive semi-definite symmetric matrix with the form
$$A=\lf(\matrix   A_{11}, & A_{12},  &\cdots, &A_{1n}\cr
          A_{21}, & A_{22},  &\cdots, &A_{2n}\cr
          \cdots, &\cdots,   &\cdots, &\cdots\cr
          A_{n1}, & A_{n2},  &\cdots, &A_{nn}
\endmatrix\ri)
$$
where $A_{ij}$ are $2\times 2$ matrix. Then
$$
\det (A)\le \Pi_{i}\det(A_{ii}).
$$
In particular, if $f$ is a real valued function defined on
$\Omega\subset \C^m$ and assumes its local minimum at a point $x\in \Omega$.
Then
$$
\det(D^2 f)(x)\le 8^m |\det(f_{\abb})|^2(x).
$$
Here $D^2f$ is the real Hessian and $f_{\abb}$ is the complex Hessian.
\endproclaim
\demo{Proof} By perturbation, one can assume that $A$ is positive definite.
Then the statement can be proved by perform the elementary row and column
operation and induction. We leave the details to the interested reader.
To prove the second part, we can choose a coordinate such that $f_{\abb}$
 is diagonal. Then
$$
|\det(f_{\abb})|^2=\lf(\frac{1}{16}\ri)^m \Pi_{i}(f_{x_ix_i}+f_{y_iy_i})^2.
$$
Since $D^2f$ is positive semi-definite, we have that
$$
(f_{x_ix_i}+f_{y_iy_i})^2\ge (f^2_{x_i x_i}+f^2_{y_i y_i})\ge
2(f_{x_ix_i}f_{y_iy_i}-f^2_{x_iy_i}).
$$
Now the result follows by applying  the first part of the lemma.
\enddemo

\demo{Proof of Proposition 1.1} For the sake of the completeness we start
with the  argument from [C].
Let $w_y(x)=R^2u(x)+\frac{1}{2}d^2_y(x)$, for $y\in B_o(R)$. Clearly we have
that $\sup_{B_o(2R)}w_y \le R^2+\frac{(3R)^2}{2}=\frac{11}{2}R^2$. But
in $B_o(7R)\setminus B_o(5R)$, we have $w_y\ge \frac{(4R)^2}{2}=8R^2\ge
\frac{11}{2}R^2$. Therefore, the minimum of $w_y$ achieves inside
$\overline{B_o(5R)}$. Namely,
$$\inf_{B_o(7R)}w_y(z)=w_y(x)$$
 for some
$x\in B_o(5R)$. As in [C] we know that at such point $x$ we have that
$$
y=\text{exp}_x(\nabla(R^2u)(x)).
$$
Now define the map $\phi(z)=\text{exp}_z(\nabla(R^2u)(z))$ for all
$z\in B_o(7R)$. We also define the measurable set
$E=\{x\in B_o(5R): \text{ there exists } y\in B_o(R) \text{ such that }
w_y(x)=\inf_{B_o(7R)}w_y\}.$ What we have shown is that for any $y\in B_o(R)$
there exists such $x\in E$ such that $\phi(x)=y$. Therefore we have that
$$
V_o(R)\le \int_E \text{Jac} (\phi)dv.
$$
Also we know for $x\in E$, $w_y(x)\le \frac{11}{2}R^2$. Therefore,
$u\le\frac{11}{2}$. This implies that
$E\subset \{u\le 6\}\cap B_o(5R)$. Therefore we further have
$$
V_o(R)\le\int_{\{u\le 6\}\cap B_o(5R)}\text{Jac} (\phi)dv.
$$
The rest is to estimate $\text{Jac} (\phi)$. The approximation argument
in [C] shows that we can apply Lemma 1.1 without the loss of the
generality even $\nabla (R^2u)(x)$ might not in $U_x$. Hence
$$
\split
\text{Jac}(\phi) &=\text{Jac}(\text{exp}_x)(R^2\nabla u)
|\det (D^2(R^2u+\frac{1}{2}d_y^2)|(x)\\
&\le |\det (D^2(R^2u+\frac{1}{2}d_y^2)|(x).
\endsplit \tag 1.9
$$
Here we have used the fact that $\text{exp}$ is volume decreasing if
$Ricci(M)\ge 0$. Now we applying Lemma 1.2 and Proposition 1.2 to estimate
$|\det (D^2(R^2u+\frac{1}{2}d_y^2)|(x)$ for $x\in E$.
$$
\split
|\det (D^2(R^2u+\frac{1}{2}d_y^2)|(x) & \le 8^m |\det (\lf(R^2 u+
\frac{1}{2}d_y^2\ri)_{\abb})|^2\\
&\le \frac{8^m}{\l^{2m}}|\det (\lf(R^2 u+
\frac{1}{2}d_y^2\ri)_{\abb})|^2|\det(a^{\abb})|^2\\
& \le \frac{8^m}{\l^{2m}}\lf\{\frac{1}{m}\sum_{\a \b}\lf(R^2 a^{\abb}u_{\abb}
+\frac{1}{2}a^{\abb}(d^2_y)_{\abb}\ri)\ri\}^{2m}\\
&\le \frac{8^m}{\l^{2m}}\lf\{\frac{1}{m}R^2Lu+\frac{1}{2}\L\ri\}^{2m}\\
& = \frac{(32)^m}{\l^{2m}}\lf\{\frac{1}{2m}R^2Lu+\frac{\L}{4}\ri\}^{2m}.
\endsplit \tag 1.10
$$
Noticing that for $x\in E$, $ \frac{1}{2m}R^2Lu+\frac{\L}{4}\ge 0$,
therefore combining (1.9) and (1.10) we have that
$$
\split
\int_{E}\text{Jac}(\phi) dv &\le
\frac{(32)^m}{\l^{2m}}\int_{E}\lf\{\frac{1}{2m}R^2Lu+\frac{\L}{4}\ri\}^{2m}
dv\\
&\le \frac{(32)^m}{\l^{2m}}\int_{\{u\le 6\}\cap B_o(5R)}\lf\{\lf(
\frac{1}{2m}R^2Lu+\frac{\L}{4}\ri)^{+}\ri\}^{2m}dv.
\endsplit
$$
This completes the proof of Proposition 1.1.

\enddemo

We close the section by pointing out that establishing a similar
result for the uniformly elliptic operator of nondivergent form on
manifolds with nonnegative Ricci curvature is an interesting
problem.

\subheading {\S 2 Monotonicity formulae}

Here we shall demonstrate a nonlinear analogy of the duality provided through
the Harnack inequality. In this case, the hinge is the monotonicity formulae.
There are various types of monotonicity formulae arising in different
geometric-analytical problems. We mainly focus on the one on harmonic maps
as well as the one on minimal submanifolds and its slightly metamorphose on
 the positive current, which is the content of the well-known
Bishop-Lelong Lemma in several complex variables. Let us start with the
harmonic maps. We refer [Sc] for the notations and proofs.

\proclaim{Proposition 2.1}
Let $u$ be a $W^{1,2}$-stationary map defined on $B_o(R)\subset \R^n$ into a
Riemannian manifold $N\subset \R^K$. Then for any point $x\in B_o(R)$ and $r>0$ such that
$B_x(r)\subset B_o(R)$ we have that
$$
I(x,r)=\frac{1}{r^{n-2}}\int_{B_x(r)}|Du|^2\, dx \tag 2.1
$$
is monotone increasing in $r$. Here $|Du|^2=\sum_{i\a}
\left(\frac{\p u^\a}{\p x^i}\right)^2$, with $u=(u^1, \cdots, u^K)
\in N \subset R^K$. Moreover we have
$$
\frac{\p }{\p \sigma}I(x,\sigma) =
2\int_{\p B_x(\sigma)}\sigma^{2-n}|\frac{\p u}{\p r}|^2\, dA. \tag 2.2
$$
\endproclaim

The above monotonicity plays an important role in the regularity theory of
the harmonic maps (cf. [S-U]). However it
 is not an easy matter
as the linear case indicated in the last section. Substantial work
are required.  One can refer [Si2] for
an updated proof of the $\e$-regularity theorem of Schoen-Uhlenbeck on
the minimizing maps.

On the other hand, analogous to the duality in the last section,
the monotonicity formula (2.2) also plays  the crucial role in
various Liouville type results on harmonic maps. For instance, there is a
result due to Garber, Ruijsenaars, Siler and Burns [G-R-S-B]
saying that any harmonic
maps $u:\R^n \to S^m$ with finite energy is a constant map for $n>3$,
which is an easy consequence of the above proposition. In [J], the author
proved that any harmonic map $u:\R^n\to N$ with $\lim_{x\to \infty} u(x)=
p_0$, a fixed point in $N$,  must be constant map. The proof again relies
crucially on Proposition 2.1. The interested readers can refer the above
mentioned papers for the detailed statement of the Liouville theorems and
their proofs. We are not going to pursue further on the
monotonicity of the harmonic maps here
and just want to mention that it would be
interesting to derive some sharp `global' version of Proposition 2.1, namely
to find
certain monotonicity which holds for a class of domain manifolds such as the
manifolds with nonnegative Ricci curvature. One expects that such a result
would be very useful in establishing and clarifying
 the  Liouville type results for the harmonic maps.

Now we want to turn our focus to the monotonicity formula for the minimal
submanifolds in $\R^n$ and its slight variation, the monotonicity formula for
 positive currents. This will be shown to be directly related to the LYH
inequality derived for the K\"ahler-Ricci flow in [N-T1] in the next
section.

Let $M^n$ be a minimal submanifold in $R^N$. The following monotonicity
formula is well-known.

\proclaim{ Proposition 2.2} Let $\Theta(x,r)$ be the density function
defined as
$$
\Theta(x,r)=\frac{1}{\omega_n r^n}\int_{B_x(r)}\,vol_M=
\frac{Vol(M\cap B_x(r))}{\omega_n r^n}. \tag 2.3
$$
Then for $r\le R$ such that $B_x(R)\cap \p M=\emptyset$,  $\Theta(x,r)$
is monotone increasing. Moreover, one has
$$
\frac{\p}{\p \sigma} \Theta(x,\sigma)=\int_{\p B_x(\sigma)}\frac{
|D^{\bot}r|^2}{\omega_n r^n}\, dA.\tag 2.4
$$
Here $\omega$ is the area of the unit sphere of dimension $n-1$,
$D^{\bot}r$ is  the normal part of $Dr$.
\endproclaim
As in [Si1], the above result also holds for stationary varifolds and it is
crucial in the regularity result  such as the the Allard regularity theorem
(cf. [Si1]). Since complex subvarieties in $\C^m$ are area minimizing,
the above
result holds for the complex subvarieties. In fact, there are slightly more
general

\proclaim{Bishop-Lelong Lemma} Let $T$ be a $(p,p)$ positive current in
$\C^m$.
Define
 $$
\Theta(x,r)=\frac{1}{r^{2m-2p}}\int_{B_x(r)}T\wedge\lf(\frac{1}{\pi}
\omega_{\C^m}\ri)^{m-p}.\tag 2.5
$$
Here $\omega_{\C^m}$ is the K\"ahler form of $\C^m$. Then
$$
\frac{\p}{\p r} \Theta(x,r)\ge 0. \tag 2.6
$$
\endproclaim

One can refer [G-H] (pages 390-391) for a proof of this lemma.
Analogous to the previous case we
 would like to relate the above monotonicity to a Liouville theorem
on plurisubharmonic functions on $\C^m$.

\proclaim{ Proposition 2.3}
let $u$ be a plurisubharmonic function defined on $\C^m$. If
$$
\lim_{x\to \infty}\frac{u(x)}{\log r(x)} =0 \tag 2.7
$$
$u$ must be a constant.
\endproclaim
\demo{Proof} The easiest proof is to restrict $u$ to $\C$ and apply the
2-dimensional result for the subharmonic functions. But we would like to
prove it using (2.6). The plurisubharmonicity of $u$ is equivalent to that
$\sqrt{-1}\ddbar u$ is a positive (1.1) current. Define
$$
M(x,r)=\frac{1}{A(x,r)}\int_{B_x(r)}u\, dA.
$$
Here $A(x,r)$ is the area of $\p B_x(r)$. Direct calculation (cf. [Ho],
Lemma 4.4.9)
shows
$$
r\frac{\p}{\p r}M(x,r)=\Theta(x,r).
$$
Namely
$$
\frac{\p}{\p \log r}M(x,r)=\Theta (x,r).
$$
Then (2.6) just implies that
$$
\frac{\p ^2}{\p t^2}M(x.r)\ge 0
$$
where $t=\log r$. Namely (2,6) is equivalent to the fact that $M(x,r)$ is a
convex function of $\log r$. Hence, under the assumption (2.7),
 we know that
$M(x,r)$ must be constant. It then implies that $\D u\equiv 0$ since
$$\Theta(x,r)=\frac{1}{r^{2m-2}}\int_{B_x(r)}\lf(\frac{1}{2}\D u\ri)
\lf(\frac{1}{\pi}\omega_{\C^m}\ri)^{m}\equiv 0.$$ The result then follows
from Theorem 1.2.
\enddemo

To prove a similar result for the plurisubharmonic functions on general
complete K\"ahler manifolds motivates the LYH inequality in the next section.

\subheading{\S 3 Linear trace Li-Yau-Hamilton inequality}
In this section we try extend the monotonicity formula (2.6) to the
general complete K\"ahler manifolds with nonnegative holomorphic
bisectional curvature. In fact, we shall show that a linear trace LYH
inequality for K\"ahler-Ricci flow would be a parabolic version of (2.6).
And this parabolic monotonicity formula is enough to prove Liouville theorem
for the plurisubharmonic functions in some cases.  Let us start with the trace version of
Hamilton's LYH inequality on the Ricci flow.

\proclaim{Theorem 3.1} Let $(M, g_{ij}(x,t))$ be a complete solution to the
Ricci flow
$$
\frac{\p}{\p t} g_{ij}(x,t)=-2R_{ij}(x,t) \tag 3.1
$$
with bounded sectional curvature and nonnegative curvature operator on
$M\times (0, T)$. Then
$$
\frac{\p R}{\p t}+2\nabla R \cdot V+2R_{ij}V_{i}V_{j}+\frac{R}{t}\ge 0
\tag 3.2
$$
for any vector field $V$.
\endproclaim
In K\"ahler category, Cao has the following
\proclaim{Theorem 3.2}
Let $(M,g_{\abb}(x,t))$ be a complete solution to the K\"ahler-Ricci flow
$$
\frac{\p}{\p t}g_{\abb}(x,t)=-R_{\abb}(x,t) \tag 3.3
$$
with bounded nonnegative holomorphic bisectional curvature on $M\times
[0, T).$ Then
$$
\frac{\p R}{\p t} +\nabla_{\bar{\a}}R V_\a +\nabla_\a RV_{\bar{\a}}+R_{\abb}
V_{\bar{\a}}V_{\b}+\frac{R}{t}\ge 0. \tag 3.4
$$
\endproclaim
Again we only state the trace version of Cao's result.
Later Chow-Hamilton extends Theorem 3.1 to the following more general case.

\proclaim{Theorem 3.3}
Let $(M, g(t))$ be a complete
solution to the Ricci flow (3.1) and let $h$ be a
symmetric 2-tensor satisfying
$$
\frac{\p}{\p t} h_{ij} =\D h_{ij}+2R_{piqj}h_{pq}-R_{ip}h_{jp}-R_{jp}h_{ip}.
\tag 3.5
$$
Assume that the metric initially has bounded
nonnegative curvature operator
(which is preserved under the flow). If on addition, $h_{ij}>0$ initially
$(h_{ij}(x,t))>0$ as long as the solution exists and
 for any vector
field $V$
$$
Q(x,t):= div(div(h))+R_{ij}h_{ij}+2div(h)\cdot V +h_{ij}V_i V_{j}+\frac{H}{2t}> 0.
\tag 3.6
$$
Here $H=g^{ij}h_{ij}$.
\endproclaim

\proclaim{Remarks.}
i) We should point out in general, one does require some growth conditions
on $h_{ij}(x,t)$ to apply the maximum principle on tensor to get
$h_{ij}(x,t)>0$. This seems missing in the statement of Theorem 3.3.
ii) Assume that the flow $(M,\ g(t))$
and the $h_{ij}$ are only defined on $M\times (0, T)$. Also assume that
the curvature operator of $g(t)$ and $\|h_{ij}(x,t)\|$ are uniformly bounded
on $M$ for any $t$. Then
the strict inequality in (3.6) will be replaced with an inequality.
The examples are the expanding solitons
 for the Ricci flow.
iii) A perturbation on $h_{ij}$ could allow the case
$h_{ij}(x,t)\ge 0$.
\endproclaim

In [NT1] the authors proved the K\"ahler version of Theorem 3.3.

\proclaim{Theorem 3.4} Let $(M, g(t))$ be  a complete solution to the
K\"ahler-Ricci flow (3.3) on $M\times [0,T)$.
Let $h_{\abb}(x,t)$ be a Hermitian symmetric
2-tensor satisfying
$$
\left(\frac{\p}{\p t}-\D \right)h_{\g\delbar}=
R_{\beta \bar{\a}\g\delbar}h_{\a\bbar}-
\frac{1}{2}\left(R_{\g\bar{p}}h_{p\delbar}+
R_{p\delbar}h_{\g\bar{p}}\right)\tag 3.7
$$
also on $M\times[0,T)$.
Assume that $(M, g(0))$ has bounded nonnegative holomorphic
bisectional curvature. Also assume that
 $(h_{\abb}(x,0))\ge 0$ and $\|h_{\abb}\|(x,0)$ is
uniformly bounded on $M$.
Then $(h_{\abb}(x,t))\ge 0$ and $Z(x,t)\ge 0$. Here $Z(x,t)$ is given by
$$
\split
Z& =g^{\a\bbar}g^{\g\delbar}\lf[\frac{1}{2}
\lf (\nabla_{\bbar}\nabla_\g +\nabla_\g \nabla_{\bbar}\ri)
h_{\a\delbar}+R_{\a\delbar}h_{\g\bbar}+\lf(\nabla_\g h_{\a\delbar}V_{\bbar}+
\nabla_{\bbar}h_{\a\delbar}V_{\g}\ri)+h_{\a\delbar}V_{\bbar}V_\g \ri]
\\
& \ \ +\frac{H}{t}.
\endsplit \tag 3.8
$$
If the $(M, g(t))$ and $h(x,t)$ are both ancient in the sense that they
exists on $M\times(-\infty, T)$. Then we have
$ \check{Z}\ge 0 $, where $\check{Z}=Z-\frac{H}{t}$.
\endproclaim

\proclaim{Remark}
The inequality $Z(x,t)\ge 0$ also holds if we only assume both the
metrics $g(t)$ and the symmetric tensor $h_{\abb}$  are only defined on
$M\times(0, T)$ provided that the curvature tensor
of $(M,\ g(t))$  and $\|h_{\abb}\|(x,t)$ uniformly bounded on $M$ for
each $t$.
\endproclaim
The Theorem 3.4 is an extension of Theorem 3.2 since in the case
$h_{\abb}=R_{\abb}(x,t)$, the Ricci tensor, (3.8) simplifies to (3.4).
The same for the real case. In [C-H], (3.6) was called linear trace Harnack
inequality. Since we have adapted the notion
Li-Yau-Hamilton inequality in [N-T1]
due to the reasons explained therein we call (3.8)
linear trace LYH inequality.

Before we explain how one can derive  a parabolic monotonicity formula out
of Theorem 3.4 we would like to indicated cases when $h_{\abb}$ satisfying
(3.7) does arise except the above mentioned obvious case that
$h_{\abb}=R_{\abb}$. Let $(M, g(t))$ be as in
Theorem 3.4. Let $(E, H)$ be a
holomorphic vector bundle on $M$ with Hermitian metric $H(x,t)$ deformed by
the Hermitian-Einstein flow:
$$
\frac{\p H}{\p t}H^{-1}=-\L F_H +\l I.\tag 3.9
$$
Here $\L$ means the contraction by  the K\"ahler form $\omega_t$,
$\l$ is a constant, which
is a holomorphic invariant in the case $M$ is compact, and $F_H$ is the
curvature of the metric $H$, which locally can be written as
$F^j_{i\abb}dz^{\a}\wedge d\bar{z}^{\b}e_{i}^*\otimes e_j$ with $\{e_i\}$ a
local frame for $E$. The transition rule for $H$ under the frame change
is $H^U_{i\bar{j}}=f_i^k\overline{f_j^k}H^V_{k\bar{l}}$ with
transition functions $f_i^j$ satisfying $e^U_i=f_i^je^V_j$.
Denote $\Omega_{\abb}=\sum_{i}F^i_{i\abb}
dz^{\a}\wedge d\bar{z}^{\b}$. One can associate $E$ the determinant line
bundle $(L, h)$. Since the transition functions for $L$ is
just $\det(f_i^j)$ we have that  $det(H^U)=|\det(f_i^j)|^2\det(H^V)$.
Namely $\det(H)$ is the naturally induced
 Hermitian metric on $L$. Using the formula
$$
F_H=\dbar(\p HH^{-1})
$$
one can easily see that $\Omega_{\abb}
dz^{\a}\wedge d\bar{z}^{\b}$ is the first Chern form of $(L, \det(H))$.
Namely $$\Omega_{\abb}
dz^{\a}\wedge d\bar{z}^{\b}=-\ddbar \log (\det(H)).$$
In this case, Theorem 3.4 can be applied
to obtain the gradient estimates for $\Omega=g^{\abb}\Omega_{\abb}$
if $(E,H)$ is a solution
of (3.9).

\proclaim{Theorem 3.5} Assume that $M$ is a compact K\"ahler manifold.
Let $(M,g(t))$ be a solution to  the K\"ahler-Ricci flow (3.3)
with nonnegative holomorphic
bisectional curvature, and let $(E,H)$ be a solution
to (3.9). Then $\Omega_{\abb}(x,t)$ satisfies (3.7) and if
$\Omega_{\abb}(x,0)\ge 0$ initially, then
$\Omega_{\abb}(x,t)\ge 0$ for all $t>0$.
Moreover, for any vector field $V$,
$$
\Omega_t+\nabla_{\a}\Omega V_{\bar{\a}}+\nabla_{\bar{\a}}\Omega V_{\a}
+\Omega_{\abb}V_{\bar{\a}}V_{\beta}+\frac{\Omega}{t}\ge 0. \tag 3.10
$$
Here $\Omega(x,t)=\sum g^{\abb}(x,t)\Omega_{\abb}(x,t).$
\endproclaim
\demo{Proof} Let us denote $\det(H(x,t))$ by $\eta(x,t)$. As we pointed out
that $\eta(x,t)$ is the induced Hermitian metric on $L$. Define
$$
u(x,t)=\log \frac{\eta(x,t)}{\eta(x,0)}.
$$
Taking trace on (3.9) then implies that
$$
\lf(\frac{\p}{\p t}H^{k}_{i}(x,t)\ri)\lf(H^{-1}\ri)^i_k=g^{\abb}
\lf\{\lf(H^k_i\ri)_{\abb}(H^{-1})^i_k+\lf(H_i^k\ri)_{\a}\lf[(H^{-1})^i_k\ri]_{\bb}\ri\}+\l r.
$$
Here $r$ is the rank of $E$. Hence
$$
\split
u_t&=\frac{\p }{\p t}\eta(x,t)\\
&= \lf(\frac{\p}{\p t}H^{k}_{i}(x,t)\ri)\lf(H^{-1}\ri)^i_k\\
&= g^{\abb}
\lf\{\lf(H^k_i\ri)_{\abb}(H^{-1})^i_k+\lf(H_i^k\ri)_{\a}\lf[(H^{-1})^i_k\ri]_{\bb}\ri\}+\l r.
\endsplit
$$
Calculating the $\D u$ writes
$$
\D u = g^{\abb}
\lf\{\lf(H^k_i\ri)_{\abb}(H^{-1})^i_k+\lf(H_i^k\ri)_{\a}\lf[(H^{-1})^i_k\ri]_{\bb}\ri\}+g^{\abb}\Omega_{\abb}(x,0).
$$
Combining them together we have
$$
\frac{\p u}{\p t}=\D u -g^{\abb}(x,t)\Omega_{\abb}(x,0)+\l r \tag 3.11
$$
 Locally we can write $-\Omega_{\abb}(x,0)$ as
$\p_{\a}\dbar_{\beta} \vp(x)$ for some smooth function $\vp$. Define locally
$U(x,t)=u+\vp -r\l t$. Then (3.11) implies
that
$$
U_t=\D U.
$$
Here $\D$ is the time-dependent Laplacian operator.
We also have that locally $U_{\abb}(x,t)=\Omega_{\abb}(x,t)$. Now Lemma 2.1
of [N-T1] implies that
$U_{\abb}(x,t)$  satisfies (3.7).  Therefore we have shown that
$\Omega_{\abb}$ satisfies (3.7). Now just apply
Hamilton maximum principle for tensors.
We then have $\Omega_{\abb}(x,t)\ge 0$ if
it is true initially.
Due to the fact that
one  can locally express $\Omega_{\abb}$ by
the $\ddbar U$ the similar reduction as carried out in the proof of
Theorem 2.1 still works. Namely if $h_{\abb}$ in Theorem 3.4 is given
by $\Omega_{\abb}=U_{\abb}$, then
$$
div(h)_{\a}=\nabla_{\g}U_{\a\bar{\g}}= \nabla_\a (U_t)
=\nabla_\a\Omega\ \  \text{ and}\ \
div(h)_{\delbar}=\nabla_{\bar{\a}}U_{\a\delbar}=\nabla_{\delbar}(U_t)
=\nabla_{\delbar}\Omega,
$$
and
$$
Z=\D (U_t) +R_{\a\bbar}U_{\bar{\a}\beta}+\nabla_{\bar{\a}}\Omega V_{\a}
+\nabla_{\a}\Omega V_{\bar{\a}}+U_{\a\bbar}V_{\bar{\a}}V_{\beta}+
\frac{\Omega}{t}
\ge 0 \tag 3.12
$$
for any $(1,0)$ vector field $V$. Here we have used the fact that
$\Omega=U_t=\D U$. Differentiate $\D U=U_t$ with respect to $t$ we have that
$$
\D (U_t)+R_{\a\bbar}U_{\bar{\a}\beta}=U_{tt}=\Omega_t.
$$
Then (3.12) implies (3.10).
\enddemo

Theorem 3.5  is just a LYH inequality
for the Hermitian-Einstein flow for vector bundles coupled with
the K\"ahler-Ricci flow. We should point out that on complete noncompact
manifold one needs some growth assumptions on $\Omega_{\abb}(x,0)$ in order
to show that $\Omega_{\abb}(x,t)\ge 0$ for later time $t>0$. For example,
it would be sufficient to assume that $\|\Omega_{\abb}\|(x,0)\le A$
for some constant $A>0$. In [N-T1], we proved the nonnegativity for
$\Omega_{\abb}(x,t)$ under much weaker assumption. Please refer to the paper
 for the details.
  The interested reader can write the
similar result for the noncompact case by modelling the
assumptions in [N-T1].
 The special case when $h_{\abb}$ in Theorem 3.4
is given by the complex Hessian of a global defined function only make sense
when $M$ is noncompact and it is clearly  the special case of the
noncompact version of Theorem 3.5. In the following we want to show that
this special case implies a parabolic monotonicity formula for the
plurisubharmonic functions.

Let us first collect some simple calculations from [N-T1]. Let $u(x)$ be a
plurisubharmonic function on $M$. Consider the heat equation coupled
with the K\"ahler-Ricci flow $\heat u(x,t)=0$ with $u(x,0)=u(x)$. With the
assumption that $\D u(x)\le C\exp (ar^2(x))$ for some positive
constants $a$ and $C$, one can show
that $u(x,t)$ is still plurisubharmonic. Applying (3.10)
we have the following
$$
w_t+\nabla_{\bar{\a}}wV_{\a}+\nabla_\a w V_{\bar{\a}}+
u_{\abb}V_{\bar{\a}}V_{\beta}+\frac{w}{t}\ge 0
$$
with $w=u_t$, which implies that
$$
w_t+\frac{w}{t}\ge 0. \tag 3.13
$$
Namely $(tw)_t \ge 0$. But
$$tw=\frac{\p \, \lf(u (x,t)\ri)}{\p\,  (\log t)}.$$ Therefore, (3.13) is
equivalent to the fact $u(x,t)$ is convex function of $\log t$. We would
like to call this fact or (3.13) the parabolic monotonicity formula for
the plurisubharmonic functions. The reason we have pointwise monotonicity
instead of  integral one is that
 the heat equation has the averaging
effect.  This more or less  says
 that the average quantity, $M(x,r)$ in the last section, is equivalent
to $u(x,t)$. For the linear heat equation on complete Riemannian manifolds
with  nonnegative Ricci curvature, the author have
proved a precise statement reflecting this principle for manifolds with
 nonnegative Ricci curvature (
cf. [N], section 3, The moment type estimates).
However due to the nonlinearity of the K\"ahler-Ricci flow
this can only be proved for the time-dependent Laplacian-heat equation
under the assumptions such as  that there exists $C>0$ with
$\aint_{B_x(r)}{\Cal R}(x,0)\, dy\le \frac{C}{r^2},$ for any $x\in M$,
 which also ensures
the long time existence of the K\"ahler-Ricci flow by the work of
W.-X. Shi [Sh2].
We also
need to assume that the bisectional curvature of  $M$ is uniformly bounded
since this is needed for the short time existence of the K\"ahler-Ricci flow
 (cf. [Sh1]).
 Under all these assumptions the same argument as in
the proof of Proposition 2.3 proves a Liouville theorem for the
plurisubharmonic functions on complete K\"ahler manifolds with nonnegative
holomorphic bisectional curvature. The interested reader can refer [N-T1] for
more detailed statements and other cases.
We should point out that using a different
method in [N-T2], the authors improved the Liouville theorems to general
complete K\"ahler manifolds with nonnegative bisectional curvature.

We would like to end this section by mentioning
 some questions arising from this
consideration on LYH inequalities. First, one might wonder what is the
`global' elliptic monotonicity? (Noting that in [N-T2] the proof  also uses
parabolic equations, can we have an elliptic proof for the Liouville
theorem on plurisubharmonic functions?) Namely what is the corresponding
 Bishop-Lelong Lemma on complete K\"ahler  manifold with nonnegative
bisectional
curvature? Certainly, the Euclidean argument does not lead to anything new.
One needs an essentially more general argument. The second question is
what is the corresponding LYH inequality for the $(p,p)$ currents. Inequality
(3.13) is the parabolic monotonicity for
the plurisubharmonic functions, or more generally  $(1,1)$
currents. What's the  parabolic monotonicity formula for
$(p,p)$-currents?
The last question is how about the real case. Namely, is there anything
analogous to what we described here for Theorem 3.3, the linear trace LYH of
Chow-Hamilton? Will it relate the minimal submanifold theory to Ricci flow
more explicitly?

\subheading{\S 4 Equalities in the  linear trace LYH inequalities}
In this section we study the geometry when
 the equality in Theorem 3.5 holds for some point at space-time. We
first prove the following result. For the simplicity of the proof we
just consider the case that
$\|h_{\abb}\|(x,t)$ are uniformly bounded in space for any fixed $t$.

\proclaim{Theorem 4.1}
Let $(M, g(t))$ be  a complete solution to the
K\"ahler-Ricci flow (3.3) on $M\times (0, T)$. Let $h_{\abb}(x,t)$ be a Hermitian symmetric
2-tensor satisfying (3.7) on $M\times (0, T)$ too.
Assume that $(M, g(0))$ has bounded nonnegative holomorphic
bisectional curvature,
 $(h_{\abb}(x,t))\ge 0$ and $\|h_{\abb}\|(x,t)$ is
uniformly bounded on $M$ for any $t\in (0, T)$. We also assumes
that $M$ is simply-connected.
Then  $Z(x,t)=0$ for some point $(x_0, t_0)$ if and only if
 $(M, g(t))$ is an expanding
K\"ahler-Ricci soliton. In the case of ancient solution, namely both flows
are  defined on $M\times (-\infty, T)$,
$\check{Z}(x,t)=0$ for some $(x_0, t_0)$ at space-time if and only if
$(M, g(t))$
is a gradient K\"ahler-Ricci soliton.
\endproclaim

This clearly implies Cao's
theorems on limits of solutions to the K\"ahler-Ricc flow in [Ca2]. The proof
is simpler, in my opinion,  and unifies the proof for
type II and type III singularity models.
The result is also more general. The similar argument  works for the
real case. Namely, by considering the equality in the LYH
inequality of
Chow-Hamilton's Theorem 3.3, we can obtain a proof of Hamilton's theorem
on the eternal (Cf. [H3]) as well as the more recent theorem on the
type III singularity
models proved by Chen-Zhu in [C-Z]. In fact we write the following
more general result.

\proclaim{Theorem 4.2}
Let $(M, g(t))$ be a complete
solution to the Ricci flow (3.1) on $M\times (0, T)$ and let $h$ be a
symmetric 2-tensor satisfying (3.5).
Assume that the metric initially has bounded
nonnegative curvature operator
(which is preserved under the flow) and $M$ is simply-connected.
We also assume that
$(h_{ij}(x,t))>0$ and $\|h_{ij}\|(x,t)$ is uniformly bounded on $M$ for any
$t\in (0, T)$. Then  $Q(x,t)=0$ for some $(x_0, t_0)$ in the
space-time if and only if
 $(M, g(t))$ is an expanding soliton. Similarly, in the case
of ancient solutions, (3.6) can be improved to
$\check{Q}\ge 0$ with $\check{Q}=Q-\frac{H}{2t}$. Moreover,
$\check{Q}(x_0, t_0)=0$ for some point $(x_0, t_0)$ in the space-time
if and only if
$(M, g(t))$ is a gradient Ricci soliton.
\endproclaim

\demo{Proof of Theorem 4.1}
First, we claim that by the maximum principle that for $t<t_0$ at every point
$x\in M$ there exists a vector $V$ such that $Z(z,t)=0$. This can be shown
using the similar argument as in Proposition 3.2 of [Ca2], where the case
for the K\"aler-Ricci flow (not the linear trace with general $h_{\abb}$) was
proved. On the other hand $Z(x,t)\ge 0$ for any arbitrary $V'$.
Therefore the $V$,
which makes $Z(x,t)=0$,  is the minimizing one. Note also here we have
assumed
$(h_{\abb}(x,t))>0$, the calculation in [N-T1], especially (1.40) shows that
$$
\split
\lf(\frac{\p}{\p t}-\D \ri) Z & = Y_1+R_{\a\bar{p}}R_{p\bar{s}}h_{s\bar{\a}}
-R_{\a\bar{p}}h_{\g\bar{\a}}\nabla_{p}V_{\bar{\g}}-R_{p\bar{\a}}h_{\a\bar{\g}}
\nabla_{\bar{p}}V_{\g}\\
&\ \ \ + h_{\g\bar{\a}}\nabla_s V_{\bar{\g}}\nabla_{\bar{s}}V_{\a}
+h_{\g\bar{\a}}\nabla_{\bar{s}}V_{\bar{\g}}\nabla_s V_\a -
\frac{H}{t^2}.
\endsplit \tag 4.1
$$
Here $Y_1$ is given by
$$
Y_1=\lf[ \D R_{s\bar{t}} +R_{s\bar{t}\a\bbar}R_{\bar{\a}\beta}
+\nabla_\a R_{s\bar{t}}V_{\bar{\a}}+\nabla_{\bar{\a}}R_{s\bar{t}}V_{\a}+
R_{s\bar{t}\a\bbar}V_{\bar{\a}}V_\beta +\frac{R_{s\bar{t}}}{t}\ri]h_{\bar{s}t}.\tag 4.2
$$
If we denote
$$
Y_2=h_{\g\bar{\a}}\lf[ \nabla_p V_{\bar{\g}}-R_{p\bar{\g}}
-\frac{1}{t}g_{p\bar{\g}}\ri]\lf[\nabla_{\bar{p}}V_\a-R_{\a\bar{p}}
-\frac{1}{t}
g_{\bar{p}\a}\ri]+
 h_{\g\bar{\a}}\nabla_{\bar{p}}V_{\bar{\g}}\nabla_pV_{\a}
$$
we can write (4.1) as
$$
\lf(\frac{\p}{\p t}-\D \ri) Z  = Y_1+Y_2-\frac{1}{t}
\lf[- h_{\g\bar{\a}}\nabla_\a V_{\bar{\g}}- h_{\g\bar{\a}}\nabla_{\bar{\a}}V_{\g}
+2R_{\a\bar{\g}}h_{\g\bar{\a}}+\frac{2H}{t}\ri].
$$
On the other hand, since $V$ now is the minimizing vector, the direct
calculation, using the equalities from the first variation consideration,
shows
$$
2Z=- h_{\g\bar{\a}}\nabla_\a V_{\bar{\g}}-
h_{\g\bar{\a}}\nabla_{\bar{\a}}V_{\g}
+2R_{\a\bar{\g}}h_{\g\bar{\a}}+\frac{2H}{t}.\tag 4.3
$$
One can refer (1.36), (1.37) and (1.39) of [N-T1]
for the detailed calculations.
Combining (4.1)--(4.3), we have
$$
\heat Z =Y_1+Y_2-\frac{2Z}{t}.
$$
On the other hand, for the minimizing $V$, from $Z(x,t)=0$ we know that
$$
\heat Z(x,t) =0.
$$
By the nonnegativity of $Y_1$ and $Y_2$ we have that
$$
Y_1=Y_2=0.
$$
Noticing that $h_{\abb}$ is positive definite, $Y_1=0$ implies that
$$
\D R_{s\bar{t}} +R_{s\bar{t}\a\bbar}R_{\bar{\a}\beta}
+\nabla_\a R_{s\bar{t}}V_{\bar{\a}}+\nabla_{\bar{\a}}R_{s\bar{t}}V_{\a}+
R_{s\bar{t}\a\bbar}V_{\bar{\a}}V_\beta +\frac{R_{s\bar{t}}}{t}=0 \tag 4.4
$$
and $Y_2=0$ implies  that
$$
\nabla_p V_{\bar{\g}}-R_{p\bar{\g}}
-\frac{1}{t}g_{p\bar{\g}}= \nabla_{\bar{p}}V_\a-R_{\a\bar{p}}
-\frac{1}{t}
g_{\bar{p}\a}  =0, \tag 4.5
$$
as well as
$$
\nabla_{\bar{p}}V_{\bar{\g}}=\nabla_pV_{\a}=0. \tag 4.6
$$
By the simply-connectness of $M$, (4.5) and (4.6) imply
that  $V$ is a holomorphic vector field and is the gradient of a
holomorphic function. Namely we have proved that $(M,\ g(t))$ is an expanding
K\"ahler-Ricci soliton.

Now suppose that $(M, g(t))$ is an expanding K\"ahler-Ricci soliton. Then we
have (4.5) above.  On the other hand, since for the minimizing vector $V$,
(4.3) holds. Plugging (4.5) into (4.3) we have
 $Z(x,t)=0$ for the minimizing vector.

The case for the ancient solution follows similarly. In this case
replace $t$ by $t-A$ and let $A\to -\infty $ we can have that
$$
\split
\lf(\frac{\p}{\p t}-\D \ri) \check{Z} & =
\check{Y}_1+R_{\a\bar{p}}R_{p\bar{s}}h_{s\bar{\a}}
-R_{\a\bar{p}}h_{\g\bar{\a}}\nabla_{p}V_{\bar{\g}}-R_{p\bar{\a}}h_{\a\bar{\g}}
\nabla_{\bar{p}}V_{\g}\\
&\ \ \ + h_{\g\bar{\a}}\nabla_s V_{\bar{\g}}\nabla_{\bar{s}}V_{\a}
+h_{\g\bar{\a}}\nabla_{\bar{s}}V_{\bar{\g}}\nabla_s V_\a
\endsplit \tag 4.1'
$$
where
$$
\check{Y}_1=\lf[ \D R_{s\bar{t}} +R_{s\bar{t}\a\bbar}R_{\bar{\a}\beta}
+\nabla_\a R_{s\bar{t}}V_{\bar{\a}}+\nabla_{\bar{\a}}R_{s\bar{t}}V_{\a}+
R_{s\bar{t}\a\bbar}V_{\bar{\a}}V_\beta \ri]h_{\bar{s}t}.\tag 4.2'
$$
Then we have that
$$
\heat \check{Z}=\check{Y}_1+\check{Y}_2
$$
with
$$
\check{Y}_2=h_{\g\bar{\a}}\lf[ \nabla_p V_{\bar{\g}}-R_{p\bar{\g}}
\ri]\lf[\nabla_{\bar{p}}V_\a-R_{\a\bar{p}}\ri]+
 h_{\g\bar{\a}}\nabla_{\bar{p}}V_{\bar{\g}}\nabla_pV_{\a}.
$$
Now the proof follows verbatim.
\enddemo

\demo{Proof of Theorem 4.2}
Similarly, using the maximum principle we know that for $t<t_0$,
at every point $x\in M$ there exists a (an unique) minimizing
vector $V$ such that $Q(x,t)=0$. One can refer Proposition 4.2 of [C-Z]
for a proof. Then the same argument as in the above proof works. One
only need to use (6.5) of Chow-Hamilton's [C-H] to replace (4.1) above.
Namely, the proof of
Chow-Hamilton's linear trace LYH (Theorem 3.3 in the last section)
 already implied
both Hamilton's Main Theorem in [H3] as well as the later Theorem 4.3 of
Chen-Zhu on the type III singularity models.
\enddemo

\medskip

{\it Acknowledgment.} This paper is based on a lecture given by the author
at the Workshop on Geometric Evolution Equations at National Center for
Theoretical Sciences at Hsinchu, Taiwan during July 15 to August 14, 2002.
The author  would like to thank the
organizers, especially Ben Chow for inviting him to the workshop as well as
 his interest to this work. The author would
also like to thank Leon Simon for a
conversation on the LYH inequality and the Liouville theorem for
plurisubharmonic functions, Jiaping Wang for bringing  Cabr\'e's result to
his attention and encouraging him to study the K\"ahler case.


\Refs
\widestnumber \key{\bf M-S-Y-1}

\ref\key{\bf B-G}\by E. Bombieri and E. Giusti
 \paper Harnack's inequality for elliptic differential equations on minimal
surfaces
\jour Invent. Math. \vol 15 \yr 1972 \pages 24--46
\endref

\ref\key{\bf C}\by X. Cabr\'e \paper Nondivergent elliptic equation on
manifolds with nonnegative curvature
\jour Comm. Pure and Appl. Math. \vol L \yr 1997 \pages623--665
\endref

\ref\key{\bf Ca1}\by H.-D. Cao \paper On Harnack inequalities for the
K\"ahler-Ricci flow
\jour Invent. Math. \vol 109 \yr 1992 \pages247--263
\endref

\ref\key{\bf Ca2}\by H.-D. Cao \paper Limits of solutions to the
K\"ahler-Ricci flow
\jour J. Differential Geom. \vol 45 \yr 1997 \pages 257--272
\endref


\ref\key{\bf C-C1}\by B. Chow and S. Chu \paper A geometric interpretation of
Hamilton's Harnack inequality for the Ricci flow
\jour Math. Res. Letter \vol 2 \yr 1995 \pages 701--718
\endref

\ref\key{\bf C-C2}\by B. Chow and S. Chu \paper A geometric interpretation to the linear
trace Harnack inequality for the Ricci flow
\jour Math. Res. Letter \vol 3 \yr 1996 \pages 549--568
\endref

\ref\key{\bf C-H}\by B. Chow and R. Hamilton\paper
Constrained and linear Harnack inequalities for
parabolic equations\jour Invent. Math. \vol 129
 \yr 1997 \pages no. 3, 213--238\endref




\ref\key{\bf C-N}\by H.-D. Cao and L. Ni \paper
Matrix Li-Yau-Hamilton estimates for heat equation on K\"ahler manifolds
\paperinfo submitted
\endref

\ref\key{\bf C-Y}\by S. Y. Cheng and S.-T.   Yau\paper Differential equations on Riemannian manifolds and their
geometric applications\jour Comm. Pure Appl. Math. \vol 28\yr 1975\pages 333--354
\endref

\ref\key{\bf C-Z}\by B.-L Chen and X. Zhu\paper
Complete Riemannian manifolds with pointwise pinched curvature
\jour Invent. Math. \vol 140
 \yr 2000 \pages 423--452\endref



\ref\key{\bf G-H}\by P. Griffith and J. Harris  \paper Principles of
Algebraic Geometry
 \paperinfo
John Wiley and Sons, 1978
\endref

\ref\key{\bf Gr}\by A. Grigor'yan \paper The heat equation on noncompact
Riemannian manifolds \jour Math. USSR Sbornik \vol 72\yr 1992\pages 47--77
\endref


\ref\key{\bf G-R-S-B}\by W.-D. Garber, S. Ruijsenaars, E. Seiler and
D. Burns \paper On finite action solution of the nonlinear $\sigma$-model
 \jour Ann. Phys. \vol 119\yr 1979\pages 305--325
\endref

\ref\key{\bf G-T}\by D. Gilbarg and N. Trudinger \paper Elliptic Partial
Differential Equations of Second Order \paperinfo 2nd Edition,.
Springer-Verlag, Berlin-Heidelberg-New York, 1983
\endref


\ref\key{\bf H1}\by R. Hamilton\paper The Harnack estimate for the Ricci
flow
\jour J. Differential Geom.\vol 37\yr 1993\pages 225--243
\endref


\ref\key{\bf H2}\by R. Hamilton \paper A matrix Harnack estimate for the
heat equation\jour Comm. Anal. Geom. \vol 1 \yr 1993 \pages 113--126 \endref

\ref\key{\bf H3}\by R. Hamilton \paper Eternal solutions to  the
Ricci flow\jour  J. Differential Geom.\vol 38 \yr 1993 \pages 1--11 \endref


\ref\key{\bf Ho}\by L. H\"ormander \paper Introduction to Complex
Analysis in Several Variables \paperinfo 3rd Edition., North-Holland, 1990
\endref


\ref\key{\bf Hu}\by G. Huisken \paper Asymptotic behavior for singularities
of the mean curvature flow
\jour J. Differential. Geom.  \vol 31 \yr 1990 \pages 285--299 \endref

\ref\key{\bf J}\by Z. Jin \paper Liouville theorems for harmonic maps
\jour Invent. Math.  \vol 108 \yr 1992 \pages 1--10 \endref



\ref\key{\bf L-W}\by P. Li and J.-P. Wang \paper
 H\"older estimates and regularity of holomorphic and harmonic functions
\jour J. Differential. Geom. \vol 58
                               \yr 2001 \pages 309-329\endref


\ref\key{\bf L-Y} \by P. Li and S.-T. Yau\paper On the parabolic kernel of
the Schr\"odinger operator\jour Acta Math.\vol 156\yr 1986\pages 139--168
\endref




\ref\key{\bf N}\by L. Ni\paper Poisson equation and Hermitian-Einstein
metrics on holomorphic vector bundles over complete noncompact K\"ahler
manifolds \jour Indiana Univ. Math. Jour. \vol51\yr 2002\pages 679--704
\endref


\ref\key{\bf N-T1} \by L. Ni and L.-F.Tam\paper Plurisubharmonic functions
and K\"ahler Ricci flow \paperinfo submitted\endref


\ref\key{\bf N-T2} \by L. Ni and L.-F.Tam\paper Louville properties of
 plurisubharmonic functions\paperinfo submitted\endref


\ref\key{\bf Sa}\by L. Saloff-Coste\paper Uniformly elliptic
operators on Riemannian manifolds \jour J. Differential Geom.\vol 36 \yr 1992
\pages 417--450 \endref

\ref\key{\bf Sc}\by R. Schoen\paper Analytic aspects of harmonic map
problem \paperinfo Seminar on nonlinear PDEs, 321--358, 1984
 \endref

\ref\key{\bf S-U}\by R. Schoen and K. Uhlenbeck\paper
A regularity theory for harmonic maps \jour  J. Differential Geom.\vol 17
 \yr 1982
\pages 307--336
 \endref

\ref\key{\bf Sh1}\by W.-X. Shi\paper Deforming the metric on complete
Riemannian manifolds \jour J. Differential Geom.\vol 30 \yr 1989
\pages 223--301 \endref

\ref\key{\bf Sh2}\by W.-X. Shi\paper Ricci deformation of metric on
complete noncompact K\"ahler manifolds \paperinfo Ph. D. thesis
Harvard University, 1990
\endref

\ref\key{\bf Si1}\by L. Simon\paper Lectures on Geometric Measure Theory
 \paperinfo Proceedings of the Centre for Mathematical Analysis, Australian
National University, {\bf 3}(1983)
\endref

\ref\key{\bf Si2}\by L. Simon\paper Theorems on Regularity and Singularity
 of Energy Minimizing Maps \paperinfo ETH
Lectures in Mathematics, Z\"urich, Birkh\"auser, 1996
\endref

\ref\key{\bf St}\by M. Struwe\paper On the evolution of harmonic maps in higher
 dimension \jour J. Differential Geom.\vol 28 \yr 1988
\pages 485--502
\endref


\ref\key{\bf  Y}\by   S.-T.   Yau\paper Harmonic functions on complete
 Riemannian manifolds\jour Comm. Pure Appl. Math.
\vol 28\yr 1975\pages 201--228
\endref


\endRefs
\enddocument